\documentclass[11pt]{amsart}
\usepackage{amsmath,amssymb,pstricks}

\newcommand{\Z}{\mathbb{Z}}
\newcommand{\R}{\mathbb{R}}
\newcommand{\C}{\mathbb{C}}

\newtheorem{theorem}{Theorem}

\newtheorem{propn}[theorem]{Proposition}
\newtheorem{conjecture}[theorem]{Conjecture}

\def\define{\textit}

\psset{unit=1pt}
\psset{arrowsize=4pt 1}
\psset{linewidth=.5pt}

\countdef\x=23
\countdef\y=24
\countdef\z=25
\countdef\t=26

\def\graybox(#1,#2){
\x=#1 \y=#2 
\z=\x \t=\y
\advance\z by 10 
\advance\t by 10 
\psframe[fillstyle=solid,fillcolor=lightgray,linewidth=0pt](\x,\y)(\z,\t) 
\psline[linewidth=.5pt](\x,\y)(\x,\t)(\z,\t)(\z,\y)(\x,\y)}

\begin{document}

\title[On quantum products of Schubert classes in a Grassmannian]{A note on quantum products of Schubert classes in a Grassmannian}
\author{Dave Anderson}
\date{July 19, 2006}
\begin{abstract}
Given two Schubert classes $\sigma_{\lambda}$ and $\sigma_{\mu}$ in the quantum cohomology of a Grassmannian, we construct a partition $\nu$, depending on $\lambda$ and $\mu$, such that $\sigma_{\nu}$ appears with coefficient $1$ in the lowest (or highest) degree part of the quantum product $\sigma_{\lambda}\star\sigma_{\mu}$.  To do this, we show that for any two partitions $\lambda$ and $\mu$, contained in a $k\times(n-k)$ rectangle and such that the $180^{\circ}$-rotation of one does not overlap the other, there is a third partition $\nu$, also contained in the rectangle, such that the Littlewood-Richardson number $c_{\lambda\mu}^{\nu}$ is $1$.
\end{abstract}
\maketitle

The purpose of this note is to establish the following fact about the product of classes in the quantum cohomology of a Grassmann manifold:
\begin{propn}\label{main}
If $d$ is the smallest or largest power of $q$ appearing in the quantum product $\sigma_{\lambda}\star\sigma_{\mu}$, then there exists a Schubert class $\sigma_{\nu}$ such that the Gromov-Witten invariant $c_{\lambda\mu}^{\nu}(d)$ is equal to $1$.
\end{propn}
In fact, we will explicitly construct such a class.  The main idea is to use a result of Postnikov (Corollary 8.4 in \cite{post}), which equates these Gromov-Witten invariants to certain classical Littlewood-Richardson numbers.  The above proposition then follows from a statement about classical cohomology (Proposition \ref{sliding} below), which says that whenever $\lambda$ and $\mu$ are such that $\sigma_{\lambda}\cdot\sigma_{\mu}\neq 0$, one can construct a partition $\nu$ such that $c_{\lambda\mu}^{\nu}=1$.  Moreover, we conjecture that the result holds for \emph{all} powers of $q$ appearing in $\sigma_{\lambda}\star\sigma_{\mu}$.  We conclude with a comment on an application of this fact to ``real quantum cohomology.''

Before discussing Postnikov's result and the construction of the class $\sigma_{\nu}$, we recall some basic definitions, notation, and results related to quantum cohomology of Grassmannans.  Let $X = Gr_k(\C^n)$ be the Grassmannian of $k$-planes in $\C^n$.  The cohomology ring $H^*(X;\Z)$ is well-understood.  It has a linear basis of \define{Schubert classes} $\sigma_{\lambda}$, indexed by partitions whose Young diagrams fit inside the $k$-by-$(n-k)$ rectangle; these classes correspond to the \define{Schubert varieties} $\Omega_{\lambda}$ of codimension $|\lambda| = \#(\text{boxes in }\lambda)$.  (Thus the class $\sigma_{\lambda}$ has degree $2|\lambda|$.)  The structure constants for multiplication in this basis are the \define{Littlewood-Richardson numbers} $c_{\lambda\mu}^{\nu}$ -- that is, 
\[ \sigma_{\lambda}\cdot\sigma_{\mu} = \sum_{\nu} c_{\lambda\mu}^{\nu}\, \sigma_{\nu}. \]

The \define{(small) quantum cohomology ring} $QH^*(X)$ is a module over the polynomial ring $\Z{[q]}$, where $q$ is a formal variable of degree $n$, with a corresponding $\Z{[q]}$-basis of Schubert classes $\sigma_{\lambda}$.  The ring structure is given by \define{quantum multiplication}, denoted by `$\star$', which has for structure constants the \define{(three-point, genus $0$) Gromov-Witten invariants}.  That is, 
\[ \sigma_{\lambda}\star\sigma_{\mu} = \sum_d q^d \sum_{\nu} c_{\lambda\mu}^{\nu}(d)\, \sigma_{\nu} ,\]
where $c_{\lambda\mu}^{\nu}(d)$ is, by definition, the number of degree-$d$ rational curves passing through general translates of $\Omega_{\lambda}$, $\Omega_{\mu}$, and $\Omega_{\nu^{\vee}}$; by degree considerations, it is nonzero only when $|\lambda|+|\mu| = |\nu| + dn$.

The ring $QH^*(X)$ has been much-studied in recent years; we mention only a few results most relevant to this note.  Agnihotri showed that the quantum product $\sigma_{\lambda}\star\sigma_{\mu}$ is never zero (see \cite{bcf}, \S5); Fulton and Woodward gave a characterization of the lowest power of $q$ appearing in $\sigma_{\lambda}\star\sigma_{\mu}$, and generalized this to all $G/P$ \cite{fw}; Yong gave an upper bound for the powers of $q$ appearing in a quantum product and conjectured that these powers form an unbroken sequence from lowest to highest; Postnikov refined the results of \cite{fw} for type $A$, gave a formula for equating the Gromov-Witten invariants $c_{\lambda\mu}^{\nu}(d)$ to Littlewood-Richardson numbers when $d$ is the minimal or maximal power of $q$ appearing in $\sigma_{\lambda}\star\sigma_{\mu}$, and proved Yong's conjecture \cite{post}.

\bigskip
Now we introduce some notation, following \cite{post}.  All partitions will lie inside the $k$-by-$(n-k)$ rectangle.  If we draw the diagram of a partition $\lambda$ inside the rectangle, the border traces a path from the $SW$ corner to the $NE$ corner of the rectangle; the \define{01-word} $\omega(\lambda)$ is the $n$-digit string which assigns a ``$0$'' to each step right, and a ``$1$'' to each step up.  Writing $\omega(\lambda) = (\omega_1,\ldots,\omega_n)$, define a doubly infinite integer sequence $\phi = \phi(\lambda) = (\phi_i)_{i\in\Z}$ by $\phi_i = \omega_1 + \cdots + \omega_i$ for $1\leq i \leq n$, and $\phi_{i+n} = \phi_i + k$ for all $i$.  Also, let $\lambda^{\vee}$ denote the \define{complement} of $\lambda$ -- that is, $\lambda^{\vee} = (n-k-\lambda_k,\ldots,n-k-\lambda_1)$ -- and let $\lambda'$ be the \define{conjugate} of $\lambda$.  Here is an example, for $k=5$ and $n=11$:

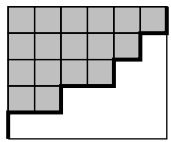
\begin{figure}[ht]
\pspicture(-20,-5)(330,60)
\psline{-}(0,0)(0,50)(60,50)(60,0)(0,0)
\graybox(0,40)
\graybox(10,40)
\graybox(20,40)
\graybox(30,40)
\graybox(40,40)
\graybox(50,40)
\graybox(0,30)
\graybox(10,30)
\graybox(20,30)
\graybox(30,30)
\graybox(40,30)
\graybox(0,20)
\graybox(10,20)
\graybox(20,20)
\graybox(30,20)
\graybox(0,10)
\graybox(10,10)
\psline[linewidth=1.5pt]{-}(0,0)(0,10)(20,10)(20,20)(40,20)(40,30)(50,30)(50,40)(60,40)(60,50)
\rput[l](80,55){$\lambda = (6,5,4,2)$}
\rput[l](80,40){$\lambda^{\vee} = (6,4,2,1)$}
\rput[l](80,25){$\lambda' = (4,4,3,3,2,1)$}
\rput[l](80,10){$\omega(\lambda) = (1,0,0,1,0,0,1,0,1,0,1)$}
\rput[l](80,-5){$\phi(\lambda) = (\ldots,1,1,1,2,2,2,3,3,4,4,5,\ldots)$}
\endpspicture
\end{figure}

Define the \define{cyclic rotation} $S^i(\lambda)$ to be the partition whose $01$-word is obtained from $\omega(\lambda)$ by cyclically permuting $i$ places to the left (or $-i$ places to the right, if $i$ is negative).  For instance, with $\lambda=(6,5,4,2)$, we have $\omega(\lambda)=(1,0,0,1,0,0,1,0,1,0,1)$, so $\omega(S^2(\lambda)) = (0,1,0,0,1,0,1,0,1,1,0)$, and thus $S^2(\lambda)=(5,5,4,3,1)$.  
Finally, given two partitions $\lambda$ and $\mu$, define integers $D_{min}$ and $D_{max}$ by
\begin{eqnarray*}
D_{min} &=& -\min_i \{ \phi_i(\lambda) + \phi_{-i}(\mu) \} \\
D_{max} &=& -\max_i \{ \phi_{-i}(\lambda) + \phi_{i - (n-k)}(\mu) \} .
\end{eqnarray*}
Of course, it suffices to consider $1\leq i \leq n$ in this definition, since the sequences $\{ \phi_i(\lambda) + \phi_{-i}(\mu) \}$ and $\{ \phi_{-i}(\lambda) + \phi_{i - (n-k)}(\mu) \}$ are $n$-periodic.

\begin{figure}[ht]
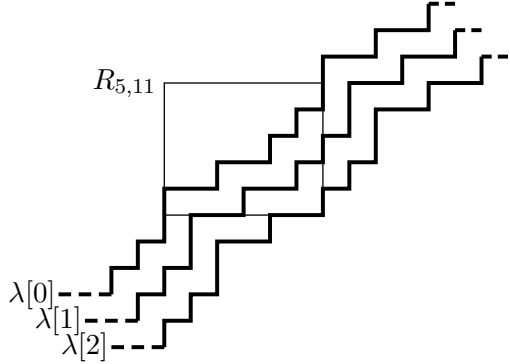

\pspicture(-140,-50)(220,80)
\rput(-15,50){$R_{5,11}$}
\psline{-}(0,0)(0,50)(60,50)(60,0)(0,0)
\rput(-50,-30){$\lambda{[0]}$}
\psline[linewidth=1.5pt,linestyle=dashed]{-}(-40,-30)(-20,-30)
\psline[linewidth=1.5pt]{-}(-20,-30)(-20,-20)(-10,-20)(-10,-10)(0,-10)(0,0)(0,10)(20,10)(20,20)(40,20)(40,30)(50,30)(50,40)(60,40)(60,50)(60,60)(80,60)(80,70)(100,70)(100,80)
\psline[linewidth=1.5pt,linestyle=dashed]{-}(100,80)(110,80)
\rput(-40,-40){$\lambda{[1]}$}
\psline[linewidth=1.5pt,linestyle=dashed]{-}(-30,-40)(-10,-40)
\psline[linewidth=1.5pt]{-}(-10,-40)(-10,-30)(0,-30)(0,-20)(10,-20)(10,-10)(10,0)(30,0)(30,10)(50,10)(50,20)(60,20)(60,30)(70,30)(70,40)(70,50)(90,50)(90,60)(110,60)(110,70)
\psline[linewidth=1.5pt,linestyle=dashed]{-}(110,70)(120,70)
\rput(-30,-50){$\lambda{[2]}$}
\psline[linewidth=1.5pt,linestyle=dashed]{-}(-20,-50)(0,-50)
\psline[linewidth=1.5pt]{-}(0,-50)(0,-40)(10,-40)(10,-30)(20,-30)(20,-20)(20,-10)(40,-10)(40,0)(60,0)(60,10)(70,10)(70,20)(80,20)(80,30)(80,40)(100,40)(100,50)(120,50)(120,60)
\psline[linewidth=1.5pt,linestyle=dashed]{-}(120,60)(130,60)
\endpspicture
\caption{Cylindric loops, for $\lambda = (6,5,4,2)$ \label{fig:loops}}
\end{figure}

The meaning of these definitions becomes clearer in the language of Postnikov's \define{toric shapes}.  (We will not need these notions for the proof of Proposition \ref{main}, but we will use them to formulate our conjecture for intermediate powers of $q$.)  Consider the lattice $\Z^2$ in the plane, with matrix coordinates; i.e., the point $(i,j)$ is $i$ steps down and $j$ steps right from the origin.  Let $R_{kn}$ be the rectangle with vertices $(0,0)$, $(k,0)$, $(0,n-k)$, and $(k,n-k)$, and let the cylinder $C_{kn}$ be the quotient $\Z^2/\Z\cdot(-k,n-k)$.  (Thus the SW and NE corners of $R_{kn}$ are identified in $C_{kn}$.)  If $\lambda$ is a partition inside $R_{kn}$, the \define{cylindric loop} $\lambda{[0]}$ is the image of the border of $\lambda$ in $C_{kn}$.  The \define{shifted cylindric loop} $\lambda{[d]}$ is the translation of $\lambda{[0]}$ by $(d,d)$.  We will often identify $\lambda{[d]}$ with its preimage in the plane; this is just the periodic continuation of the (translated) border of $\lambda$.  See Figure \ref{fig:loops}.

\begin{figure}[ht]
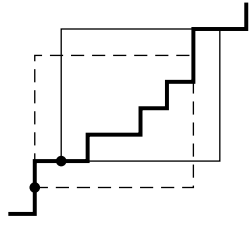

\pspicture(-150,0)(200,65)
\psline[linestyle=dashed]{-}(0,0)(0,50)(60,50)(60,0)(0,0)
\psline[linewidth=1.5pt]{-}(-10,-10)(0,-10)(0,0)(0,10)(20,10)(20,20)(40,20)(40,30)(50,30)(50,40)(60,40)(60,50)(60,60)(80,60)(80,70)
\psline[linestyle=solid]{-}(10,10)(10,60)(70,60)(70,10)(10,10)
\pscircle*(0,0){2}
\pscircle*(10,10){2}
\endpspicture
\caption{Rotating the frame \label{fig:frames}}
\end{figure}

A \define{frame} is any translation of $R_{kn}$ in the plane, and the \define{anchor} of a frame is its SW corner.  If we move a frame so that its anchor lies on $\lambda{[0]}$, then the part of $\lambda{[0]}$ contained inside the frame forms the border of a partition.  In fact, if the anchor is shifted $i$ steps in the NE direction along $\lambda{[0]}$, then the resulting partition is $S^i(\lambda)$.  Also, the number $\phi_i(\lambda)$ is the vertical distance traveled after $i$ steps NE along $\lambda{[0]}$ (so the frame for $S^i(\lambda)$ is translated up by $\phi_i(\lambda)$ from $R_{kn}$, and right by $i - \phi_i(\lambda)$).  See Figure \ref{fig:frames}.

If $\lambda$ and $\mu$ are partitions such that $\mu{[d]}$ is (weakly) right and below $\lambda{[0]}$ in the plane, so that the region between $\mu{[d]}$ and $\lambda{[0]}$ forms a connected strip, then the image of this region in $C_{kn}$ is called a \define{cylindric shape} and denoted $\mu/d/\lambda$.  Let $\lambda^{\downarrow}{[0^{\downarrow}]}$ denote the translation of $\lambda{[0]}$ by $(k,0)$.  A cylindric shape $\mu/d/\lambda$ is \define{toric} if $\mu{[d]}$ lies between $\lambda{[0]}$ and $\lambda^{\downarrow}{[0^{\downarrow}]}$.  It is not hard to see that the numbers $D_{min}$ and $D_{max}$ defined above are the minimum and maximum values of $d$ such that $\mu^{\vee}/d/\lambda$ is a toric shape.

\begin{figure}[ht]
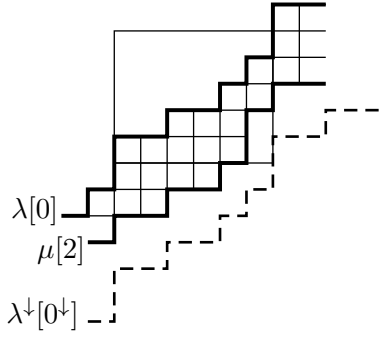

\pspicture(-140,-55)(220,70)
\psline{-}(0,0)(0,50)(60,50)(60,0)(0,0)
\rput(-30,-18){$\lambda{[0]}$}
\psline[linewidth=1.5pt]{-}(-20,-20)(-10,-20)(-10,-10)(0,-10)(0,0)(0,10)(20,10)(20,20)(40,20)(40,30)(50,30)(50,40)(60,40)(60,50)(60,60)(80,60)
\rput(-20,-33){$\mu{[2]}$}
\psline[linewidth=1.5pt]{-}(-10,-30)(0,-30)(0,-20)(20,-20)(20,-10)(40,-10)(40,0)(50,0)(50,20)(60,20)(60,30)(80,30)
\rput(-27,-57){$\lambda^{\downarrow}{[0^{\downarrow}]}$}
\psline[linewidth=1pt,linestyle=dashed]{-}(-10,-60)(0,-60)(0,-50)(0,-40)(20,-40)(20,-30)(40,-30)(40,-20)(50,-20)(50,-10)(60,-10)(60,0)(60,10)(80,10)(80,20)(100,20)
\psline{-}(0,-20)(0,0)
\psline{-}(10,-20)(10,10)
\psline{-}(20,-10)(20,10)
\psline{-}(30,-10)(30,20)
\psline{-}(40,0)(40,20)
\psline{-}(50,20)(50,30)
\psline{-}(60,30)(60,40)
\psline{-}(70,30)(70,60)
\psline{-}(-10,-20)(0,-20)
\psline{-}(0,-10)(20,-10)
\psline{-}(0,0)(40,0)
\psline{-}(20,10)(50,10)
\psline{-}(30,20)(50,20)
\psline{-}(40,30)(60,30)
\psline{-}(60,40)(80,40)
\psline{-}(60,50)(80,50)
\psline{-}(60,60)(80,60)
\endpspicture
\caption{The toric shape $\mu/2/\lambda$, for $\mu=(4,3,3,2)$ and $\lambda=(6,5,4,2)$. \label{fig:toricshape}}
\end{figure}

Postnikov shows that $q^d$ appears in the quantum product $\sigma_{\lambda}\star\sigma_{\mu}$ if and only if $\mu^{\vee}/d/\lambda$ is a toric shape, and deduces that $D_{min}$ and $D_{max}$ are, respectively, the minimum and maximum powers of $q$ appearing in this product.  By the definitions, there are integers $a$ and $b$ such that $D_{min} + \phi_a(\lambda) + \phi_{-a}(\mu) = 0$ and $D_{max} + \phi_{-b}(\lambda) + \phi_{b-(n-k)}(\mu) = 0$.  (There may be more than one such $a$ and $b$, but any choice will do.)  Set 
\begin{eqnarray*}
\lambda^{min} &=& S^a(\lambda), \\
\mu^{min} &=& S^{-a}(\mu), \\
\lambda^{max} &=& S^b(\lambda^{\vee}), \quad\text{ and} \\
\mu^{max} &=& S^{n-k-b}(\mu^{\vee}) .
\end{eqnarray*}
Then Postnikov proves the following:
\begin{propn}[{\cite{post}, Corollary 8.4}] \label{p}
Let $\nu$ be any partition in the $k$-by-$(n-k)$ rectangle.  Then 
\begin{eqnarray}
c_{\lambda\mu}^{\nu}(D_{min}) &=& c_{\lambda^{min}\mu^{min}}^{\nu},\text{ and}\\
c_{\lambda\mu}^{\nu}(D_{max}) &=& c_{\lambda^{max}\mu^{max}}^{\nu^{\vee}}.
\end{eqnarray}
\end{propn}
In particular, the products $\sigma_{\lambda^{min}}\cdot\sigma_{\mu^{min}}$ and $\sigma_{\lambda^{max}}\cdot\sigma_{\mu^{max}}$ are nonzero in $H^*(X)$.

By substituting $\lambda$ for $\lambda^{min}$, and so on, this reduces Propostion \ref{main} to the following:
\begin{propn} \label{sliding}
Let $\lambda$ and $\mu$ be any partitions contained in the $k$-by-$(n-k)$ rectangle, such that $\sigma_{\lambda}\cdot\sigma_{\mu} \neq 0$ in $H^*(X)$.  Then there is a partition $\nu=\nu(\lambda,\mu)$, also contained in the rectangle, such that $c_{\lambda\mu}^{\nu} = 1$.
\end{propn}
If we write $\mu_{180}$ for the $180^{\circ}$-rotation of $\mu$ inside the $k$-by-$(n-k)$ rectangle, note that the condition $\sigma_{\lambda}\cdot\sigma_{\mu} \neq 0$ is equivalent to requiring that $\lambda$ and $\mu_{180}$ do not overlap.  (This notation should cause no confusion, as we will not discuss partitions with $180$ parts.)  Note that the boxes of $\mu_{180}$ form the complement of $\mu^{\vee}$ inside the rectangle.
\begin{proof}
We will construct the partition $\nu$, and use the following version of the Littlewood-Richardson rule: The number $c_{\lambda\mu}^{\nu}$ is equal to the number of semistandard Young tableaux of shape $\nu/\lambda$ with reading word of type $\mu$.\footnote{The \define{reading word} of a tableau is the integer string formed by reading the entries of the tableau from right to left, starting at the top row.  A word $w=w_1 w_2 \cdots w_p$ is of type $\mu$ if one can build the diagram of $\mu$ by placing a box in row $w_1$, then in row $w_2$, etc., in such a way that one has a Young diagram at each step.  The condition that each stage be a Young diagram is equivalent to requiring that for each $m\leq p$, 
\[ \#(1\text{'s in }\{w_1,\ldots,w_m\}) \geq \#(2\text{'s in }\{w_1,\ldots,w_m\}) \geq \cdots;\] this is called the \define{Yamanouchi condition}.}  We will call such a tableau on $\nu/\lambda$ a \define{Littlewood-Richardson filling} of type $\mu$.  (See \cite{yt} or \cite[Appendix 1]{ec2} for more on the Littlewood-Richardson rule.)

Draw $\lambda$ and $\mu_{180}$ inside the rectangle.  Now slide the columns of $\mu_{180}$ up against $\lambda$, and then left-justify all rows.  The resulting shape is $\nu(\lambda,\mu)$.  Here is an example, with $k=5$, $n=11$, $\lambda = (4,3,1)$, and $\mu = (5,4,4)$.  (The shape of $\lambda$ is shaded, and that of $\mu_{180}$ is filled by numbers.)

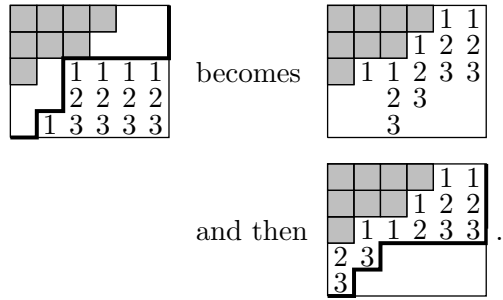
\pspicture(-60,-65)(280,70)
\psline{-}(0,0)(0,50)(60,50)(60,0)(0,0)
\graybox(0,20)
\graybox(0,30)
\graybox(10,30)
\graybox(20,30)
\graybox(0,40)
\graybox(10,40)
\graybox(20,40)
\graybox(30,40)
\rput(15,5){1}
\rput(25,5){3}
\rput(35,5){3}
\rput(45,5){3}
\rput(55,5){3}
\rput(25,15){2}
\rput(35,15){2}
\rput(45,15){2}
\rput(55,15){2}
\rput(25,25){1}
\rput(35,25){1}
\rput(45,25){1}
\rput(55,25){1}
\psline[linewidth=1.5pt]{-}(0,0)(10,0)(10,10)(20,10)(20,30)(60,30)(60,50)

\rput*[l](70,25){becomes}

\psline{-}(120,0)(120,50)(180,50)(180,0)(120,0)
\graybox(120,20)
\graybox(120,30)
\graybox(130,30)
\graybox(140,30)
\graybox(120,40)
\graybox(130,40)
\graybox(140,40)
\graybox(150,40)
\rput(135,25){1}
\rput(145,5){3}
\rput(155,15){3}
\rput(165,25){3}
\rput(175,25){3}
\rput(145,15){2}
\rput(155,25){2}
\rput(165,35){2}
\rput(175,35){2}
\rput(145,25){1}
\rput(155,35){1}
\rput(165,45){1}
\rput(175,45){1}

\rput*[l](70,-35){and then}

\psline{-}(120,-60)(120,-10)(180,-10)(180,-60)(120,-60)
\graybox(120,-40)
\graybox(120,-30)
\graybox(130,-30)
\graybox(140,-30)
\graybox(120,-20)
\graybox(130,-20)
\graybox(140,-20)
\graybox(150,-20)
\rput(135,-35){1}
\rput(125,-55){3}
\rput(135,-45){3}
\rput(165,-35){3}
\rput(175,-35){3}
\rput(125,-45){2}
\rput(155,-35){2}
\rput(165,-25){2}
\rput(175,-25){2}
\rput(145,-35){1}
\rput(155,-25){1}
\rput(165,-15){1}
\rput(175,-15){1}
\psline[linewidth=1.5pt]{-}(120,-60)(130,-60)(130,-50)(140,-50)(140,-40)(180,-40)(180,-10)

\rput*(185,-38){.}
\endpspicture

\noindent
In this example, then, $\nu(\lambda,\mu) = (6,6,6,2,1)$.

This sliding algorithm is reminiscent of the moves in Sch\"utzenberger's jeu de taquin \cite{schutz} (see also \cite[Appendix 1]{ec2}).  In fact, the bulk of the sliding described here can be accomplished via jeu-de-taquin moves; however, as the above example shows, it is not exactly the same as jeu de taquin.  (In jeu de taquin, the `$3$' in the bottom row would slide up, and the final shape would be $(6,6,6,3)$.)

Numerically, let $\rho$ be the partition formed by sorting 
\[ (k- \lambda'_1 - \mu'_{n-k}, k - \lambda'_2 - \mu'_{n-k-1}, \ldots, k - \lambda'_{n-k} - \mu'_1) . \]
(These are the heights of the spaces between the columns of $\lambda$ and $\mu_{180}$.)  The sliding construction described above leaves the shape $(\rho')_{180}$ in the bottom right corner.  Indeed, sliding $\mu_{180}$ up leaves blank columns of heights $(k - \lambda'_i - \mu'_{n-k+1-i})$, and left-justifying the filled space is the same as right-justifying blank space, which is equivalent to sorting.
Thus $\nu(\lambda,\mu) = (\rho')^{\vee}$.  In the above example, $\rho = (2,2,2,2,1)$, so $\rho' = \nu^{\vee} = (5,4)$. 

Now we must show that $c_{\lambda\mu}^{\nu}=1$.  First, we exhibit a Littlewood-Richardson filling of $\nu/\lambda$, proving $c_{\lambda\mu}^{\nu}\geq 1$.  In fact, the tableau produced in our running example is a Littlewood-Richardson filling; we claim the procedure suggested there works in general.  Let us make this precise.  Consider $\mu_{180}$ as a skew shape, and fill its boxes by writing the numbers $1,2,3,\ldots$ down columns, so that the $r$th column from the right has entries $1,2,\ldots,\mu'_r$.  Note that this is a Littlewood-Richardson filling of type $\mu$.  Now slide the boxes as prescribed (first moving them up against $\lambda$, then left-justifying), carrying their labels along.  The result is, by definition, a tableau on the shape $\nu/\lambda$.

We need to check that the result is actually a Littlewood-Richardson filling of type $\mu$.  By construction, the tableau has entries corresponding to $\mu$.  The sliding operations preserve weak increase along rows and strict increase down columns, so the tableau is semistandard.  It remains to verify the Yamanouchi condition; for this, we will consider the intermediate shape $\theta$ formed by sliding $\mu_{180}$ up against $\lambda$, and the corresponding filling of $\theta$ -- this is obtained by filling the columns of $\theta$ just as was done with $\mu_{180}$, so that the $r$th column from the right has entries $1,\ldots,\mu'_r$.  Note that the reading word is unchanged by left-justification, so it suffices to show that the reading word of this filling (of $\theta$) satisfies the Yamanouchi condition.

Let $B$ be the $m$th box one reads when forming the reading word $w$.  The letters $w_1,\ldots,w_m$ are the entries appearing in rows strictly above $B$, or in the same row and weakly right of $B$.  In Figure \ref{fig:word}, $B$ is the darkly shaded box, and the entries in question are all those in the shaded region.  Every entry in a given column is distinct, so the number of $i$'s apearing in the shaded region is bounded by the number of columns in the shaded region.  There is a $1$ at the top of each column, so we see that 
\[ \#(1\text{'s}) = \#(\text{columns}) \geq \#(i\text{'s}) \]
for each $i>1$.  If we remove the boxes filled with $1$'s, we can repeat this argument on the part of the shaded region that remains; this shows that the Yamanouchi condition holds.

\begin{figure} 
\pspicture(-100,0)(200,65)
\psline{-}(0,0)(0,60)(90,60)(90,0)(0,0)
\pspolygon[linewidth=.5pt,fillstyle=none,fillcolor=lightgray](0,15)(5,15)(5,30)(10,30)(10,35)(20,35)(20,40)(35,40)(35,45)(55,45)(55,60)(0,60)
\pspolygon[linewidth=0pt,linecolor=lightgray,fillstyle=solid,fillcolor=lightgray](15,30)(15,35)(20,35)(20,40)(35,40)(35,45)(55,45)(55,60)(90,60)(90,25)(60,25)(60,35)(55,35)(55,25)(50,25)(50,30)
\pspolygon*(50,25)(50,30)(55,30)(55,25)
\pspolygon[linewidth=1pt](15,25)(15,35)(20,35)(20,40)(35,40)(35,45)(55,45)(55,60)(90,60)(90,15)(75,15)(75,20)(70,20)(70,25)(60,25)(60,35)(55,35)(55,20)(40,20)(40,25)(30,25)(30,30)(25,30)(25,25)
\rput(72,47){$\theta$}
\rput(12,47){$\lambda$}
\endpspicture
\caption{\label{fig:word}}
\end{figure}

\medskip
One can prove the reverse inequality $c_{\lambda\mu}^{\nu}\leq 1$ by pondering tableaux, but here is a simpler way, pointed out to me by Sergey Fomin.  Let $\rho$ be the sorting of the numbers $(k-\lambda'_i-\mu'_{n-k+1-i})$, as above.   First, note that $\rho'_1$ is the size of the (unique) largest horizontal strip which can be added to $\lambda$ without overlapping $\mu_{180}$ or spilling outside the rectangle.  Indeed, $\rho'_1$ is the number of nonzero parts of $\rho$, which is the number of columns (of the rectangle) in which there is space between $\lambda$ and $\mu_{180}$.  It follows (by Pieri's rule) that $\sigma_{\lambda} \cdot \sigma_{\rho'_1} \cdot \sigma_{\mu} = \sigma_{\tilde{\lambda}} \cdot \sigma_{\mu}$, where $\tilde{\lambda}$ is the shape formed by adding this longest horizontal strip to $\lambda$.  If we write $\tilde{\rho}$ for the partition formed from the vertical spaces between $\tilde{\lambda}$ and $\mu_{180}$, then $\tilde{\rho}'_1 = \rho'_2$.  Proceeding inductively, we see that 
\begin{eqnarray}
\sigma_{\lambda} \cdot (\sigma_{\rho'_1}\cdot \cdots \cdot \sigma_{\rho'_s}) \cdot \sigma_{\mu} &=& \sigma_{\mu^{\vee}} \cdot \sigma_{\mu} \nonumber \\
&=& 1\cdot{[pt]}.
\end{eqnarray}

It follows that $c_{\lambda\mu}^{\alpha^{\vee}} \leq 1$ for every partition $\alpha$ appearing in the Schubert expansion of $(\sigma_{\rho'_1}\cdot\cdots\cdot\sigma_{\rho'_s})$.  Since $\rho'$ is such a partition, and $\nu = (\rho')^{\vee}$, we are done.
\end{proof}

\begin{figure}[ht]
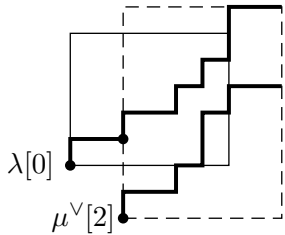

\pspicture(-140,-20)(220,55)
\psline{-}(0,0)(0,50)(60,50)(60,0)(0,0)
\rput(-15,0){$\lambda{[0]}$}
\psline[linewidth=1.5pt]{-}(0,0)(0,10)(20,10)(20,20)(40,20)(40,30)(50,30)(50,40)(60,40)(60,50)(60,60)(80,60)
\rput(5,-20){$\mu^{\vee}{[2]}$}
\psline[linewidth=1.5pt]{-}(20,-20)(20,-10)(40,-10)(40,0)(50,0)(50,20)(60,20)(60,30)(80,30)
\pscircle*(0,0){2}
\pscircle*(20,-20){2}
\pscircle*(20,10){2}
\psline[linestyle=dashed]{-}(20,-20)(20,60)(80,60)(80,-20)(20,-20)
\endpspicture
\caption{Setup for sliding, with $\lambda=(6,5,4,2)$ and $\mu=(6,4,3,3,2)$. \label{fig:conj1}}
\end{figure}

\medskip
We conclude by describing a conjectured algorithm for producing a class $\nu = \nu(\lambda,\mu,d)$, for each $d$ between $D_{min}$ and $D_{max}$, such that $c_{\lambda\mu}^{\nu}(d)=1$.  Begin by drawing the paths $\lambda{[0]}$ and $\mu^{\vee}{[d]}$; mark the point on $\mu^{\vee}{[d]}$ which is the translation of the anchor by $(d,d)$.  (See Figure \ref{fig:conj1}.)  Consider the box formed by the union of two frames: one whose anchor is at $(d,d)$, and the other whose anchor is at the point of $\lambda{[0]}$ directly above $(d,d)$.
Perform the sliding algorithm described in the proof of Proposition \ref{sliding} for the shapes whose borders are the parts of $\lambda{[0]}$ and $\mu^{\vee}{[d]}$ lying inside this box.  Call the partition produced by the sliding algorithm $\tilde{\nu}$, and let $\nu(\lambda,\mu,d)$ be the partition formed by the last $k$ parts of $\tilde{\nu}$ (including zeroes).  (This is the part of $\tilde{\nu}$ lying inside the frame whose anchor is at $(d,d)$.)

For example, with $k=5$, $n=11$, $\lambda=(6,5,4,2)$, and $\mu=(6,4,3,3,2)$, the algorithm produces $\nu(\lambda,\mu,2)=(6,6,1)$:

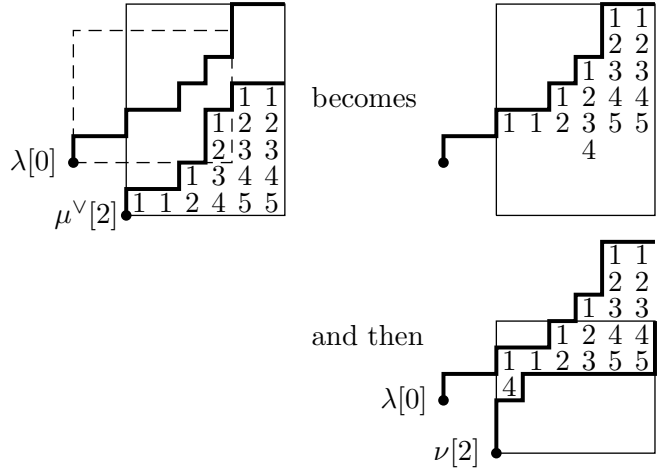
\pspicture(-60,-120)(280,120)
\psline[linestyle=dashed]{-}(0,0)(0,50)(60,50)(60,0)(0,0)
\rput(-15,0){$\lambda{[0]}$}
\psline[linewidth=1.5pt]{-}(0,0)(0,10)(20,10)(20,20)(40,20)(40,30)(50,30)(50,40)(60,40)(60,50)(60,60)(80,60)
\rput(5,-20){$\mu^{\vee}{[2]}$}
\psline[linewidth=1.5pt]{-}(20,-20)(20,-10)(40,-10)(40,0)(50,0)(50,20)(60,20)(60,30)(80,30)
\pscircle*(0,0){2}
\pscircle*(20,-20){2}
\psline[linestyle=solid]{-}(20,-20)(20,60)(80,60)(80,-20)(20,-20)
\rput(25,-15){1}
\rput(35,-15){1}
\rput(45,-15){2}
\rput(55,-15){4}
\rput(65,-15){5}
\rput(75,-15){5}
\rput(45,-5){1}
\rput(55,-5){3}
\rput(65,-5){4}
\rput(75,-5){4}
\rput(55,5){2}
\rput(65,5){3}
\rput(75,5){3}
\rput(55,15){1}
\rput(65,15){2}
\rput(75,15){2}
\rput(65,25){1}
\rput(75,25){1}

\rput*[l](90,25){becomes}

\psline[linewidth=1.5pt]{-}(140,0)(140,10)(160,10)(160,20)(180,20)(180,30)(190,30)(190,40)(200,40)(200,50)(200,60)(220,60)
\pscircle*(140,0){2}
\psline[linestyle=solid]{-}(160,-20)(160,60)(220,60)(220,-20)(160,-20)
\rput(165,15){1}
\rput(175,15){1}
\rput(185,15){2}
\rput(195,5){4}
\rput(205,15){5}
\rput(215,15){5}
\rput(185,25){1}
\rput(195,15){3}
\rput(205,25){4}
\rput(215,25){4}
\rput(195,25){2}
\rput(205,35){3}
\rput(215,35){3}
\rput(195,35){1}
\rput(205,45){2}
\rput(215,45){2}
\rput(205,55){1}
\rput(215,55){1}

\rput*[l](90,-65){and then}

\rput(125,-90){$\lambda{[0]}$}
\pscircle*(140,-90){2}
\psline[linewidth=1.5pt]{-}(140,-90)(140,-80)(160,-80)(160,-70)(180,-70)(180,-60)(190,-60)(190,-50)(200,-50)(200,-40)(200,-30)(220,-30)
\pscircle*(160,-110){2}
\psline[linestyle=solid]{-}(160,-110)(160,-60)(220,-60)(220,-110)(160,-110)
\rput(165,-75){1}
\rput(175,-75){1}
\rput(185,-75){2}
\rput(165,-85){4}
\rput(205,-75){5}
\rput(215,-75){5}
\rput(185,-65){1}
\rput(195,-75){3}
\rput(205,-65){4}
\rput(215,-65){4}
\rput(195,-65){2}
\rput(205,-55){3}
\rput(215,-55){3}
\rput(195,-55){1}
\rput(205,-45){2}
\rput(215,-45){2}
\rput(205,-35){1}
\rput(215,-35){1}
\rput(145,-110){$\nu{[2]}$}
\psline[linestyle=solid, linewidth=1.5pt]{-}(160,-110)(160,-90)(170,-90)(170,-80)(220,-80)(220,-60)

\rput(225,-68){.}
\endpspicture

\noindent
One can check that $c_{(6,5,4,2),(6,4,3,3,2)}^{(6,6,1)}(2) = 1$.  We conjecture that this always works: if $\nu=\nu(\lambda,\mu,d)$ is as described above, for $D_{min}\leq d\leq D_{max}$, then $c_{\lambda\mu}^{\nu}(d)=1$.  In particular, we expect the following generalization of Proposition \ref{main} to hold:

\begin{conjecture} \label{conj}
If $d$ is any power of $q$ appearing in the quantum product $\sigma_{\lambda}\star\sigma_{\mu}$, then there exists a Schubert class $\sigma_{\nu}$ such that the Gromov-Witten invariant $c_{\lambda\mu}^{\nu}(d)$ is equal to $1$.
\end{conjecture}

\noindent

\medskip
To summarize, we have seen that 
\[ c_{\lambda\mu}^{\nu(\lambda^{min},\mu^{min})}(D_{min}) = c_{\lambda^{min},\mu^{min}}^{\nu(\lambda^{min},\mu^{min})} = 1\]
and 
\[ c_{\lambda\mu}^{\nu(\lambda^{max},\mu^{max})^{\vee}}(D_{max}) = c_{\lambda^{max},\mu^{max}}^{\nu(\lambda^{max},\mu^{max})} = 1 .\]
Of course, this implies that the mod $2$ reduction of $\sigma_{\lambda}\star\sigma_{\mu}$ is always nonzero.  This can be seen as an analogue of one of the main results of \cite{fw} for ``mod $2$ real quantum Schubert calculus,'' at least for Grassmannians.\footnote{The phrase in quotes should be interpreted as follows: Let $\overline{M}=\overline{M}_{0,3}(X,d)$ be the Kontsevich moduli space of stable maps, and let $\overline{M}(\R)$ be its real part.  The Gromov-Witten invariants $c_{\lambda\mu}^{\nu}(d)$ are certain intersection numbers in $H^*(\overline{M},\Z)$; let $\overline{c}_{\lambda\mu}^{\nu}(d)$ be the analogous intersection numbers in $H^*(\overline{M}(\R),\Z/2\Z)$.  It is reasonable to expect that $\overline{c}_{\lambda\mu}^{\nu}(d) \equiv c_{\lambda\mu}^{\nu}(d) \pmod 2$, as is true for the classical case ($d=0$).  An outline discussion of intersection theory on $\overline{M}(\R)$ can be found in \cite{kwon}.} 
Similarly, a proof of Conjecture \ref{conj} would establish a real analogue of the stronger result that the powers of $q$ appearing in a quantum product form an unbroken sequence from $D_{min}$ to $D_{max}$ \cite[Theorem 8.1]{post}.

\bigskip
I would like to thank William Fulton for suggesting this question and for comments on the manuscript, and Sergey Fomin for a helpful discussion.  Anders Buch's Littlewood-Richardson calculator\footnote{Available at \texttt{http://www.math.rutgers.edu/\~{}asbuch/lrcalc/}.} proved invaluable for experimentation.

\small{\textsc{Department of Mathematics, University of Michigan, Ann Arbor, MI 
48109}

\textit{E-mail address}: \texttt{dandersn@umich.edu}}

\end{document}